\documentclass[graybox]{svmult}

\usepackage{type1cm}        
\usepackage{makeidx}         
\usepackage{graphicx}        
\usepackage{multicol}        
\usepackage[bottom]{footmisc}

\usepackage{newtxtext}       %
\usepackage[varvw]{newtxmath}       

\usepackage{hyperref}

\usepackage[sort]{natbib}

\makeindex             

\begin{document}

\title*{Recent Analytical and Computational Developments on the  Advection–Diffusion-Reaction Wildfire Model}
\titlerunning{Recent Developments on the ADR Wildfire Model}
\author{Luca~Nieding\orcidID{0009-0001-8321-2281}, A.~George~Morgan\orcidID{1111-2222-3333-4444}, Adrian~Navas\orcidID{0000-0002-3465-6898}, Donato~Pera\orcidID{0000-0003-4872-4713}, Bruno~Rubino \orcidID{0000-0001-6546-8390}, Federica~Di~Michele\orcidID{0000-0003-1244-6942}, Koondanibha~Mitra\orcidID{0000-0002-8264-5982}, Qiyao~Peng\orcidID{0000-0002-7077-0727}, and \\ Cordula~Reisch\orcidID{0000-0003-1442-1474} }
\authorrunning{Nieding et al.}

\institute{
Luca Nieding \at Institute for Partial Differential Equations, TU Braunschweig, Braunschweig, Germany. \email{l.nieding@tu-braunschweig.de}
\and A. George Morgan \at York University, Department of Earth and Space Science and Engineering, Toronto, ON, Canada. \email{agmorgan@my.yorku.ca}
\and Adrian Navas \at Fluid Dynamics Technologies Group, Aragon Institute of Engineering Research (I3A), University of Zaragoza, Spain. \email{anavas@unizar.es}
\and Donato Pera, Bruno Rubino \at Department of Information Engineering
Computer Science and Mathematics
University of L'Aquila, L'Aquila, Italy. \email{donato.pera@univaq.it, bruno.rubino@univaq.it}, 
\and
Federica Di Michele \at Istituto Nazionale di Geofisica e Vulcanologia - Milano, Italy
\email{federica.dimichele@ingv.it}
\and Koondanibha Mitra (Correspondence) \at CASA, Department of Mathematics and Computer Science, Eindhoven University of Technology, the Netherlands. \email{k.mitra@tue.nl}
\and Qiyao Peng (Correspondence) \at (1) MARS, School of Mathematical Sciences, Lancaster University, Lancaster, United Kingdom; (2) Mathematical Institute, Faculty of Science, Leiden University, Leiden, The Netherlands. \email{qiyao.peng@lancaster.ac.uk}
\and Cordula Reisch (Correspondence) \at Dep. of Science and Environment, Roskilde University, Roskilde, Denmark. \email{cordular@ruc.dk}}
%
%
\maketitle

\vspace{-3em}
\abstract{Wildfires represent a problem for ecosystems, human activities, and economies, driven by the climate crisis and land-use changes. Predicting wildfire propagation through mathematical modelling is essential for damage mitigation and risk assessment. This paper provides a comprehensive review of a physics-based Advection-Diffusion-Reaction (ADR) model, focusing on the balance between physical accuracy and computational efficiency. We analyze the  ability of the ADR model to estimate fire front speed and behaviour and discuss its preliminary mathematical properties. Additionally, we discuss some modelling improvements which enhance the physical realism of the model.  Furthermore, we address the challenge of reducing computational costs, emphasizing the need for inexpensive but precise
numerical schemes. We report recent findings outlining open challenges in model discretization and technological solutions. All these developments highlight the potential of ADR models as powerful tools for efficient  wildfire simulation and risk assessment.}

\vspace{-1em}
\section{Introduction}
\vspace{-1em}
\subsection{Background and Motivation}
\vspace{-1em}

Wildfires, as a natural disaster, can cause widespread damage to the environment, ecosystems, human residential areas, and local economies. Due to the climate crisis and land-use changes, the frequency and severity of wildfires has significantly increased in the past few decades \cite{owid-wildfires,Cunningham2024}. Every year, hundreds of hectares of woodland and green areas are destroyed by fires. Although natural fires are not necessarily a problem, as they can be part of the natural cycle of an ecosystem, the repeated and excessive number of fires, whether caused by natural causes, human activity or arson, can pose a serious threat to the survival of an ecosystem and cause pollution. 

The spread of wildfires is influenced by many factors, e.g. moisture, wind, landscape, and vegetation type. Mathematical modelling may take those factors into account and predict the propagation of the fire. A feasible mathematical model for wildfire requires an \textit{almost} accurate mathematical description of the reality, and a computationally efficient numerical algorithm. Analysing the model provides insight into various outcomes depending on the parameters, particularly it brings the opportunity to estimate the speed of the travelling wave (e.g. fire front). Alongside this, with a computationally efficient model, one can obtain predictive results and take timely actions before the wildfire spreads further causing additional damages.

Various types of models, based on different methodologies, have been proposed over the years, such as data-driven, mechanistic surrogate and physics-based models, ranging from stochastic to deterministic \cite{grasso,perry1998current}.  In this review, we limit ourselves to a specific advection-diffusion-reaction (ADR) derived in \cite{Asensio2002}, and analysed in detail e.g. in \cite{REISCH2024179, grasso, sero2002modelling} and the references mentioned in the following subsections. We will summarize the most recent findings regarding this ADR model of wildfires from both mathematical and computational perspectives, and the open challenges in the field.

\vspace{-1.5em}
\subsection{The standard ADR model}
\vspace{-1em}
The typical ADR model equations are based on the conservation laws of mass and energy applied to a control volume describing the vegetation layer, under a set of simplifying assumptions \cite{sero2002modelling, navas2025modelling}. A continuum approach is used within the vegetation layer to describe the fuel-air mixture (which is assumed to be incompressible). This allows for the definition of bulk parameters that characterize its effective properties. Furthermore, a local thermal equilibrium between the fuel and the surrounding gas is assumed. Vertical energy fluxes due to volatilization are not explicitly taken into account \cite{navas2025modelling}. Thus, the ADR model equations in differential form are given by
\begin{equation}
    \label{Eq_wildfire_model}
    \left\{
    \begin{aligned}
        \rho c(T_t+\boldsymbol{v}\cdot\nabla T) &= \nabla\cdot[K(T)\nabla T]-h(T-T_\infty) + \rho\Psi_TSY,\\
        Y_t &= -\Psi_TY, 
    \end{aligned}
    \right.
\end{equation}
with the bulk temperature $T$ and the biomass fraction $Y>0$. The model parameters are the bulk density $\rho$, the specific heat at constant pressure $c$, both of which are assumed constant in most models; the advection speed $\boldsymbol{v}$ due to wind, the heat exchange coefficient with the environment $h$, the ambient air temperature $T_\infty$ and the heating value $S$. The first term on the right describes both the effect of heat diffusion and radiation, with the resulting heat transfer coefficient 
\begin{equation*}
    K(T)= k + 4\epsilon\delta\sigma T^3, 
\end{equation*}
with $k$ the heat conduction coefficient,  $\epsilon$ the emissivity factor, $\delta$ the optical path length and $\sigma$ the Stefan-Boltzmann constant, see \cite{Asensio2002} for details.
The second term models the heat loss to the environment, approximated by Newton's cooling law, where $h>0$ is the heat exchange coefficient with the environment. The last term represents the combustion process, which is described by the nonlinear function 
\begin{equation}\label{eq:Psi_1st}
 \Psi_T(x,t)= A H(T(x,t)-\bar{T})  \exp\left(-\frac{T_{\rm ac}}{T(x,t)}\right),
 \end{equation}
 with the Heaviside function $H:\mathbb{R}\to\{0,1\}$. The function $\Psi_T$ is zero for $T<\bar{T}$ and one for $T\geq \bar{T}$, where the activation temperature $T_{\rm ac}$, the ignition temperature $\Bar{T}$, and a pre-exponential factor $A>0$ arise from the Arrhenius formula.  
 \par System \eqref{Eq_wildfire_model} is completed with the initial datum
 \begin{align}\label{eq:ini}
         &T|_{t=0} = T_0(x),  \quad  Y|_{t=0}= Y_0(x), &&\quad \forall \ x\in \mathbb{R}^d.
 \end{align}
In the wildfire context, they satisfy $T_0(x)\geq T_\infty$ and $Y_0(x)>0$ for all $x\in \mathbb{R}^d$. 

\vspace{-1em}
\section{Mathematical Analysis}
\vspace{-1em}
In this section, we focus on the theoretical analysis on System \eqref{Eq_wildfire_model}-\eqref{eq:Psi_1st}, regarding the existence and uniqueness of the solution and travelling wave (e.g. the wildfire front). For the sake of analysis, further simplifications and approximations are made accordingly. 
\vspace{-2em}
\subsection{Existence-Uniqueness Theory}
\vspace{-1em}
Naturally, much of the analytical work on reaction-diffusion models of wildfire spread is devoted to travelling combustion waves (see the discussion below). However, understanding how reaction-diffusion models describe the evolution of ``arbitrary" initial states is also important. Such investigations allow us to see if there are interesting classes of solutions beyond travelling waves, and in particular if the model admits solutions with finite-time blowup; see \cite[Ch. 3]{BE1989} for an example of a combustion model with such singular solutions. 
\par Well-posedness results for wildfire models are available only in the regime when $\Psi_T$ remains Lipschitz continuous with respect to $T$ (owing to issues related to identification of the $\Psi_T$ term with its limit when applying a compactness argument). Due to the extremum principle, this is guaranteed when $\bar{T}\leq T_\infty$, since then the temperature remains over $T_\infty$ and thus $H(T-\bar{T})=1$. The well-posedness of nonlinear (possibly degenerate) PDE-ODE models proven in \cite{mitra2025quasilinear} for bounded domains is stated as Theorem \ref{thm:well-posedness}:
\begin{theorem}[Mitra \& Sonner \cite{mitra2025quasilinear}]\label{thm:well-posedness} 
Let $\bar{T}\leq T_\infty$. Suppose for a bounded Lipschitz domain $\Omega\subset \mathbb{R}^d$, an initial datum $\left(T_0, Y_0\right)\in L^{\infty}_{x}(\Omega)\times L^{\infty}_{x}(\Omega)$ satisfies $T_0\geq T_\infty$ and $Y_0\geq 0$ a.e. Then, a unique weak solution $(T,Y)$ of \eqref{Eq_wildfire_model} exists, where $T\in (L^2_t H^1_0(\Omega))\cap H^1_t H^{-1}(\Omega)\cap L^\infty_{t,x}$, $Y\in C^0_t L^2(\Omega)\cap L^\infty_{t,x}$, $T\geq T_\infty$ and $Y\geq 0$ a.e.
\end{theorem}

\noindent In the above, notation for the $L^p$ spaces and the Sobolev space $H^1$ have been used.

For unbounded domains, Morgan \cite{Morgan2025} developed some basic existence-uniqueness theory for a variant of \eqref{Eq_wildfire_model} introduced by Weber et al \cite{WMSG1997} characterized by
$$
\epsilon = \overline{T} =  T_\infty=0, \;\;  \boldsymbol{v} =\boldsymbol{0} \quad \text{(in appropriate units).}
$$
This means we are ignoring the influence of radiative transfer and drift due to wind. Next, after conducting the non-dimensionalization of the model, all other constants remaining in the problem are given by 
$$
\rho = c = k = A = S = T_{\text{ac}} = 1. 
$$
This reduces the Cauchy problem for \eqref{Eq_wildfire_model} to 
\begin{equation}\label{eqn:weber_model_ivp_full}
\left\{
\begin{aligned}
    T_{t} &=  \left(\Delta - h\right) T + \Psi_{T}Y, \quad \forall \ t>0, \ x\in \mathbb{R}^d,
    \\
    Y_{t} &= -\Psi_{T}Y,  \quad \forall \ t>0, \ x\in \mathbb{R}^d,
\end{aligned}
\right.
\end{equation}
with initial datum \eqref{eq:ini}. Additionally, we remark that $\Psi_{T}$ actually becomes \emph{smooth} (but not real-analytic) in this case (since $\overline{T} = T_{\text{ac}} =0$).
\par Note in particular that \emph{unbounded} spatial domain is required for travelling wave solutions to exist. If \eqref{eqn:weber_model_ivp_full} was instead posed on a bounded domain and the initial states were given by $T_0  = Y_0 \equiv 1$, one could immediately apply the results of \cite[Ch. 4]{BE1989} to obtain a (potentially non-unique) solution to \eqref{eqn:weber_model_ivp_full} valid on $t\in (0,\infty)$. However, accommodating an unbounded domain and generic initial conditions is straightforward: using a routine fixed-point argument in conjunction with $L^p$-estimates for the linear heat equation, we obtain
\begin{theorem}[Morgan \cite{Morgan2025}]\label{thm:linfty_gwp} 
Suppose we are given an initial datum $\left(T_0, Y_0\right)\in L^{\infty}_{x}\times L^{\infty}_{x}$. Then, the mild formulation of \eqref{eqn:weber_model_ivp_full}
\begin{subequations}\label{eqn:mild_soln}
\begin{align}
    T(x,t) &= e^{-h t}e^{t\Delta}T_0 + \int_{0}^{t} e^{-h\left(t-\tau\right)}e^{\left(t-\tau\right)\Delta} Y(x,\tau)\Psi_{T}\left(x,\tau\right) \ \mathrm{d} \tau \label{eqn:mild_T}
    \\
    Y(x,t) &= Y_0(x) \exp\left(-\int_0^t \Psi_{T}\left(x,\tau\right) \ \mathrm{d} \tau\right). 
    \label{eqn:mild_Y}
\end{align}
\end{subequations}
admits a unique solution $(T,Y)\in C^0_t L^{\infty}_{x}\times C^0_t L^{\infty}_{x}$ valid for at least some small time interval $[0,t^{*}]$. Additionally, if $h>0$ then this solution is valid for $t\in [0,\infty)$ and obeys the non-blowup estimate
    \begin{equation}
        \left\|T(x,t)\right\|_{L^{\infty}_{x}} \leq e^{-h t}\left\|T_0(x)\right\|_{L^{\infty}_{x}} + h^{-1}\left\|Y_0(x)\right\|_{L^{\infty}_{x}} \quad \forall \ t>0.
        \label{eqn:non_blowup}
    \end{equation}
If $h=0$, then the local solution is valid on $\left[0, t^{\Omega}\right]$ for \emph{any} $t^{\Omega}>0$, but \eqref{eqn:non_blowup} must be replaced with the weaker estimate 
\begin{equation}
     \left\|T(x,t)\right\|_{L^{\infty}_{x}} \leq\left\|T_0(x)\right\|_{L^{\infty}_{x}} + t\left\|Y_0(x)\right\|_{L^{\infty}_{x}}, \quad \forall \ t\in \left[0, t^{\Omega}\right].
\end{equation}
\end{theorem}
\noindent In particular, finite-time thermal blowup is impossible regardless of whether the environment-vegetation heat exchange coefficient $h$ vanishes. Note that mild solutions $\left(T, Y\right)$ are smooth in $x$ and continuously differentiable in $t$ for $t>0$, so they are almost classical solutions to \eqref{eqn:weber_model_ivp_full} (up to technicalities associated with the convergence to $\left(T_0, Y_0\right)$ as $t\rightarrow 0^{+}$). For more details, we refer to \cite{Morgan2025}.

\vspace{-2em}
\subsection{Travelling Waves}\label{sec:TW}
\vspace{-1em}
In line with real-world observations of wildfires, numerical simulations of the model in \eqref{Eq_wildfire_model} show travelling wave solutions. The formal mathematical proof is nevertheless challenging because of the nonlinear coupling of the system variables, including an unsteady step-function. 
The existence of travelling wave solutions was proven for simplified models, e.g. for \eqref{eqn:weber_model_ivp_full} in \cite{WMSG1997}, or without a step-function in \cite{VV2002}, where the stability of the wave front was aslo analysed. A pure combustion diffusion system was studied in \cite{Ghazaryan&c2010}.

An almost full model, compare \eqref{Eq_wildfire_model}, was proven in \cite{BBH2009}, where only the cooling by Newton's law was not included ($h=0$). 
In this simplified model, the main mechanism for temperatures to decrease after the burning is missing. Therefore, the temperature increases due to the combustion and then remains at its maximal value. 
In \cite{REISCH2024179}, \eqref{Eq_wildfire_model} was studied under the assumption of a constant biomass fraction $Y$ using a linearization approach. Numerical simulations compare the estimated travelling wave speed using the linearization approach for a constant biomass fraction, the simplified model in \cite{BBH2009}, and numerical simulations of the model \eqref{Eq_wildfire_model} for various parameters. The main result is that the upper and lower limits gained through the simplified model in \cite{BBH2009} provide good estimates for the numerical spreading speed as long as the variations of the parameters are small enough. 
Further, the linearization approach approximates the dependency of the rate of spread on the Newton cooling parameter $h$ better, resulting for fixed $Y=1$ in a close upper bound. 

In \cite{mitra_studying_2025}, the authors modified the combustion function to
\begin{align}\label{eq:new_psi2}
     \Psi_T(x,t)&=A_L {H\left(\Theta_{T}(x,t)-\bar{T}\right)} (T-T_\infty). 
 \end{align}
This describes a process where a fire keeps burning even if the temperature decreases, and where the Arrhenius law can be approximated by a first-order polynomial. The existence of travelling wave solutions with a reaction rate given by \eqref{eq:new_psi2} has been proved. Furthermore, with the travelling wave analysis and introducing the new energy term, the PDE system can be converted to an ODE system. A shooting algorithm is presented in \cite{mitra_studying_2025} to estimate the travelling wave speed. With this algorithm, numerical simulations show good agreement between the ordinary differential equation model and the advection diffusion reaction model.

\vspace{-1em}
\section{Model Extensions}
\vspace{-1em}
To better model and investigate how other factors (e.g. topography and environmental conditions) influence the wildfire progression, variations in the model \eqref{Eq_wildfire_model} have been investigated. During the study of the system, numerical investigations showed that the Arrhenius function in the combustion term, \eqref{eq:Psi_1st}, may be simplified. 
For analytical purposes, it was approximated in \cite{mitra_studying_2025} by a first-order polynomial, see Sec. \ref{sec:TW}, showing a good agreement in the travelling wave speed but an overestimation of the temperature. The replacement with a constant factor was studied in \cite{navas2025modelling}, showing a good agreement of both, the maximal temperature and the travelling wave speed. 

\vspace{-2em}
\subsection{Two-phase approach for bulk parameters}
\vspace{-1em}
Wildfire impact prediction is challenging because fire propagation involves multi-scale physical and chemical processes. Predictive models must explicitly resolve only the larger scales while approximating smaller-scale phenomena to remain computationally feasible and faster than real time. Following the criteria of Sero et al. \cite{sero2002modelling}, four scales are distinguished: gigascopic, macroscopic, mesoscopic, and microscopic. Fire propagation is herein modelled at the macroscopic scale, while mesoscopic and smaller-scale heterogeneities and processes are incorporated indirectly through averaged properties, allowing the treatment of the vegetation layer as a continuous porous medium with a spatial discretization on the order of the mesoscopic scale.

The wildfire propagation model is formulated over a three-dimensional domain representing the ground surface and a finite height above it. Using a continuum approach, the biomass and surrounding gas (air and combustion products) can be modelled as a two-phase porous medium composed of solid fuel and a gaseous phase. By averaging the quantities over the vertical direction, the problem is reduced to two spatial dimensions defined on the ground surface. At each point in space, bulk properties such as density, specific heat, volume, and mass are defined for both phases using mixture rules \cite{sero2002modelling,navas2025modelling}. 

When using the two-phase porous flow-based formulation, the bulk velocity $\mathbf{v}$ is related to the actual wind velocity,  $\mathbf{w}$, as  
\begin{equation}\label{eq:advectTerm} 
\mathbf{v}=  (1-R_f)\frac{\bar{\rho}_a \bar{c}_{p,a}}{\bar{\rho}_f \bar{c}_{p,f} R_{f} +  \bar{\rho}_a \bar{c}_{p,a} (1-R_f)}\mathbf{w},
\end{equation}
with $R_f$ the fuel volume fraction, and $\bar{\rho}_a$, $\bar{\rho}_f$, $\bar{c}_{p,a}$, $\bar{c}_{p,f}$ the density and specific heats of the gas and solid phases.
This bulk velocity represents the physical relations within the fuel and replaces the general advection velocity $\mathbf{v}$ in \eqref{Eq_wildfire_model}.

\vspace{-2em}
\subsection{Modelling wind and topography through advection}
\vspace{-1em}
Both wind and topography can have a strong influence on a wildfire's spread, direction, and thus also on its danger \cite{BBH2009, math8101674}. Both these effects into account using an advective term in the governing PDE system.
Now, in \eqref{Eq_wildfire_model}, the natural wind speed is already included by the advective term involving $\boldsymbol{v}$. It was shown in \cite{BBH2009,mitra_studying_2025} that for relatively small advection velocities there are two travelling wave solutions: a fast wave and a slow wave. However, at higher wind speeds, the slower wave is suppressed, and the fast wave is dominant, moving in the same direction as the wind.
Furthermore, it has been shown that the speed of the fire front increases with steeper terrain inclines \cite{sanchez-monroy_fire_2019, abouali_analysis_2021}.
  As a simple model, we describe this direction-dependent acceleration by a virtual wind $\nabla Z(\boldsymbol{x})$  where $Z(\boldsymbol{x})$ defines the topography \cite{math8101674}.
  Thus, the new advection velocity is given by
\begin{equation}
    \boldsymbol{v} = \beta \boldsymbol{w} + \gamma \nabla Z(\boldsymbol{x}),
\label{eq:virtualWind}
\end{equation}
where $\boldsymbol{w}$ is the natural advection velocity in the area and $\beta$, $\gamma$ are calibration factors (note that $\beta$ can be derived from \ref{eq:advectTerm})  Numerical simulations show that this model extension qualitatively replicates the directional dependence of the fire front on topography \cite{math8101674,NIEDING2025103}. If using a conservative version of the energy equation, the non-zero divergence of $\boldsymbol{v}$ may produce unphysical solutions. One way to solve this problem is to add a third component to the wind: a vertical top-flow. This top-flow can be chosen to cancel the non-zero divergence. Alternatively, if solving the non-conservative version of the equation, this problem is avoided \cite{NIEDING2025103}. Although this model extension is phenomenologically rather than physically motivated, it nevertheless captures the physical effect of topography on wildfire spread in a wide parameter set, \cite{NIEDING2025103}.

\vspace{-2em}
\subsection{Modelling fuel moisture}
\vspace{-1em}
One of the most significant controls on both wildfire risk and propagation in many regions of the world is the  fuel moisture content \cite{Ruffault:2022b}, defined as the ratio between the mass of water, $m_w$, and the mass of dry fuel, $m_0$, as
\begin{equation}
    M=\frac{m_w}{m_0}.
\end{equation} 
Fuel moisture content is inversely related to the likeliness of ignition and fire intensity. Although fuel moisture plays a critical role, it is often omitted from ADR-based wildfire propagation models, reducing predictive accuracy and limiting insight into its effects. 

Only a few models explicitly include the effect of fuel moisture. The PhyFire model by Asensio et al. \cite{asensio2023historical} includes sensible and latent water heating using a multi-valued enthalpy–temperature relationship, formulated as a Stefan problem \cite{eyres1946calculation, swaminathan1993enthalpy}. This approach has been validated with observations and underpins an operational fire-spread simulator \cite{ASENSIO2023105710, prieto2017gis}. Vogiatzoglou et al. \cite{vogiatzoglou2024interpretable} explicitly model the endothermic phase through two consecutive Arrhenius reactions for wood dehydration and combustion, offering greater physical detail than standard ADR models. Margerit and Séro-Guillaume \cite{margerit2002modelling, sero2002modelling, sero2008large} propose a hierarchy of increasingly complex models that account for moisture effects using stage-dependent equations and a Dirac-type term to represent latent heat at evaporation. Finally, Navas-Montilla et al. \cite{navas2025modelling} use the apparent heat capacity method, which represents thermal phase-change effects through an effective, typically piecewise-defined, heat capacity \cite{caggiano2018reviewing, swaminathan1993enthalpy}. Three moisture models of increasing complexity are introduced in \cite{navas2025modelling} to describe evaporation and fuel composition. The most comprehensive model captures the key physicochemical mechanisms governing moisture in both live and dead fuels, enabling independent estimation of their specific heat.

The simplest version of the models in Navas-Montilla et al. \cite{navas2025modelling} introduces a fuel specific heat  defined as a piecewise constant function of the temperature by
 \begin{equation}  \label{eq:cpfsimple} 
  \bar{c}_{p,f}= \left\{ \begin{array}{ccccl}
      c_{p, \mathrm{eff}} & \hbox{if} & T<\bar{T} & \hbox{ and }  & Y=1 ,\\
      c_{p,f_0} & \hbox{if} & T\geq\bar{T}  &  \hbox{ or } & Y<1, 
     \end{array}\right.
 \end{equation} 
where $\bar{T} $ is the ignition temperature, $c_{p,f_0}$ is the specific heat of the dry fuel and 
\begin{equation}\label{cpeffsimple2}
  c_{p, \mathrm{eff}}= c_{p,f_0}  + M\left[\frac{c_w(T_w-T_\infty) + L_w }{(\bar{T}-T_\infty)}\right] 
\end{equation}
is the apparent specific heat before ignition that accounts for the presence of water, with $M$ the moisture content, $c_w$ specific heat of water, $L_w$ the latent heat of evaporation of water and $T_w$ the temperature of evaporation of water.
This fuel specific heat replaces the general $\rho$ in \eqref{Eq_wildfire_model} and may model spatially heterogeneous vegetation of different types.

\vspace{-2em}
\section{Efficient Computational Methods}\label{sec:computing}
\vspace{-1em}

The interplay of the different mechanisms in the advection-diffusion-reaction model \eqref{Eq_wildfire_model} presents challenges to efficient and precise numerical simulations. We will reflect on considerations of numerical methods and computing architectures from personal experience, resulting in (not exhaustive) guidelines.

A crucially important purpose of wildfire model simulation is to provide an assisting tool for planning firefighter operations and reducing the harm of spreading fires to humans and the ecosystem. The computation time is therefore crucial. It is influenced by the choice of the mathematical model, the numerical methods, the computational infrastructure, and their interplay.  
The interplay of methods is particularly relevant when using High-Performance Computing (HPC) techniques, which can increase the quality of numerical simulations through the use of finer spatial and temporal discretization in reduced computational time. 

We distinguish HPC architectures between CPU-based (Central Processing Unit) architectures and GPU-based (Graphics Processing Units) architectures. CPU-based architectures can be divided into two further macro categories: distributed memory architectures and shared memory architectures. Each kind of macro architecture requires specific software tools and techniques for programming.

Below, we present an HPC-oriented classification of numerical methods based on two properties: the possibility of writing the numerical scheme in an explicit semi-discrete form (EX) and the high-order accuracy with hp-adaptivity (HHP) \cite{dg}.

Based on these criteria, we can analyse the properties of the most common numerical methods, such as: finite-difference method (FDM), finite volume methods (FVM), finite element methods (FEM), and discontinuous Galerkin method (DG). 
From the analysis of the two selection  criteria (EX and HHP) for each mentioned numerical method, the FDM and DG methods are able to satisfy both properties, whereas the FVM methods have problems with the HHP property and the FEM methods with the EX property. 
We will therefore focus our analysis on the FD and DG methods. 

FD methods have a lower computational cost compared to DG methods and are easy to implement on HPC infrastructures. A drawback of FDM is the difficult implementation of complex (three-)dimensional geometries. FDM is efficient for the production of frequently used two-dimensional fire-spread maps.

DG methods guarantee better performance than FD methods when dealing with complex geometries, e.g. three-dimensional ones. However, implementation in HPC structures is more difficult, and the computational cost is higher compared to FDM. 

High-order FVM are as well a suitable tool, in particular when using essentially non-oscillatory (ENO) and weighted essentially non-oscillatory (WENO) schemes.
WENO methods rely on a nonlinear, selective stencil selection (or weighting) procedure that achieves high-order accuracy in smooth regions while robustly capturing discontinuities without generating spurious Gibbs oscillations.

WENO methods require the use of high-order integration methods in time, such as Runge-Kutta schemes or Arbitrary high-order derivative (ADER) approaches, which guarantee a globally high-order accuracy in both space and time. The combination of the WENO method with the Strong-Stability preserving Runge-Kutta 3 (SSPRK3) has been applied to the ADR wildfire model in \cite{navas2025modelling}. The essentially non-oscillatory property of the WENO method is advantageous when handling under-resolved travelling waves with sharp gradients as those present in the ADR wildfire propagation model. 
Compared to the DG reconstruction, WENO methods typically exhibit a lower algorithmic complexity, while preserving robustness and simplicity in the treatment of discontinuous solutions. This choice, however, comes at the expense of a lower formal accuracy per degree of freedom and a non-local formulation, which compromises computational efficiency, as the wider stencils required in WENO can limit parallel scalability. 

A comparison between CPU and GPU architectures for FD methods showed an advantage of GPU architectures in terms of the calculation times and computational efficiency, \cite{HPC2}.
The evaluation of the performance of the other numerical methods on HPC structures will be subject of future study.

An alternative to HPC techniques may be reduced order models, where in the case of proper orthogonal decompositions, the solution is projected to a finite basis space. Previously, the standard model \eqref{Eq_wildfire_model} has been investigated.
To cover the effects of the nonlinear terms, additional interpolation methods (DEIM) may be used, leading to speed ups of the computation time of a factor around 100, see e.g. \cite{LBGLL2016}, being even improved by parameter-dependent reduced order models, see \cite{BSU2021}.
These results were obtained for one spatial dimension and heterogeneous vegetation. 
A test case in two spatial dimensions is presented in \cite{solan-fustero_paramatrized_2025} with a speed up of a factor 45, but requiring a high offline-phase effort for an acceptable accuracy of the reduced order model solution.

\vspace{-1em}
\section{Conclusion and Discussion}
\vspace{-1em}
\par When it comes to analytical problems, there are several interesting avenues for future investigations. 
First, the discontinuity in $\Psi_T$ due to the Heaviside function is still unaccounted for, imposing a big restriction on our results. Secondly, while extending Theorem \ref{thm:linfty_gwp} for unbounded domains to include a prescribed wind velocity is straightforward, radiative effects are more difficult to accommodate since they introduce \emph{nonlinear terms involving derivatives of $T$}. Additionally, one could attempt to add the effects of spotting (see, for instance \cite{MH2016}) into the model \eqref{eqn:weber_model_ivp_full} via a nonlinear, nonlocal term.

\textbf{Open challenges:} The complexity of wildfire physics naturally brings in challenges in modelling. Atmospheric conditions, vegetation type, and topography of the land all have a significant impact on the wildfire spread. As a result, even the simplest model will contain many parameters as well as complicated terms, such as non-linearities, which lead to the difficulties analysing the characteristics of the model.

Further, modelling these effects leads to complex, nonlinear equations requiring advanced computation techniques.
Future studies will focus on evaluating the performance of the other numerical methods mentioned in section \ref{sec:computing} using HPC structures together with an urgent HPC-class data acquisition and processing model.

Data availability is another challenge in modelling wildfires because accurate predictions depend on detailed, up-to-date information that is often incomplete, sparse, or inconsistent. With many parameters involved in the model, one needs data including all the information, such as weather conditions, topographic distribution of vegetation, vegetation type and its state (e.g. fuel moisture). In addition, precise ignition location is not always recorded due to safety risks and logistical constraints. These gaps can then reduce predictive reliability of the model.

Coupling the ADR wildfire model with an atmospheric model would make it possible to represent wind–fire interactions (e.g. pyroconvection). While this integration requires substantial computational resources, it offers the potential for more reliable forecasts. It would support a deeper investigation of plume-driven fires, which are expected to increase in frequency under ongoing climate change.

\vspace{-1em}
\begin{acknowledgement}
  The work of QP was supported by Research England under the Expanding Excellence in England (E3) funding stream, which was awarded to MARS: Mathematics for AI in Real-world Systems in the School of Mathematical Sciences at Lancaster University. The work of ANM was supported by the Spanish Ministry of Science, Innovation and Universities - Agencia Estatal de Investigación (10.13039/501100011033) and FEDER-EU under project-nr. PID2022-141051NA-I00. 
This work received financial support from ICSC – Centro Nazionale di Ricerca in High Performance Computing, Big Data and Quantum Computing, funded by European Union - NextGenerationEU under the Italian Ministry of University and Research (MUR) National Centre for HPC, Big Data and Quantum Computing CN 00000013-CUP:E13C2200100000.  
\end{acknowledgement}

\vspace{-1em}
\ethics{Competing Interests}{
The authors have no conflicts of interest to declare that are relevant to the content of this chapter.}

\vspace{-2em}
\bibliography{main}
\bibliographystyle{abbrv} 

\end{document}